\documentclass[12pt,reqno]{amsart}

\usepackage{latexsym}
\usepackage{amssymb}
\usepackage{graphics}

\renewcommand{\epsilon}{\varepsilon}

\newcommand{\PP}{{\mathbb P}}

\newcommand{\G}{{\mathbb G}}

\newcommand{\C}{{\mathbb C}}
\newcommand{\Q}{{\mathbb Q}}

\newcommand{\CP}{\C\PP}

\newcommand{\rank}{{\operatorname{rank}}}

\renewcommand{\phi}{\varphi}

\newcommand{\ocal}{\mathcal{O}}

\newtheorem{theo}{{Theorem}}[section]

\newtheorem{prop}[theo]{{Proposition}}

\title{On the structure of Demailly-Semple invariant jet differentials }

\author{Jingzhou Sun }
\address{Department of Mathematics, Johns Hopkins University, Baltimore, MD
21218, USA} \email{jzsun@math.jhu.edu}

\thanks{Research  partially supported by NSF grant
 DMS-0901333.}

\date{\today}

\begin{document}

\begin{abstract}

\end{abstract}

\maketitle

In the terminology of
\cite{de}, a directed manifold is a pair $(X, V)$, where $X$ is a complex manifold
and $V\subset T_X$ a subbundle.
Let $(X, V)$ be a complex directed manifold, $
J_kV \rightarrow X$ is defined to be the bundle of $k$-jets of germs of curves $f: (\C,0)
\rightarrow X$ which are
tangent to V, i.e., such that $f'(t)\in  V_{f(t)}$ for all $t$ in a neighborhood of 0, together
with the projection map $f\rightarrow f(0)$ onto $X$. It is easy to check that $J_kV$ is actually a
subbundle of $J_kT_X$.  Let $\G_k$ be the group of germs of $k$-jet
biholomorphisms of $(\C, 0)$, that is, the group of germs of biholomorphic maps
$$t \rightarrow \phi(t)
= a_1t + a_2t^2 +\cdots +a_kt^k, \qquad a_1\in \C^*, a_j\in \C, \quad j > 2$$
in which the composition law is taken modulo terms $t^j$ of degree $j > k$. The
group $\G_k$ acts on the left on $J_kV$ by reparametrization, $(\phi,f)\rightarrow f\circ\phi$.

 Given a directed manifold $(X,V)$ with $\rank V=r$, let $\tilde{X}=\PP(V)$. The subbundle $\tilde{V}\subset T_{\tilde{X}}$ is defined by $$\tilde{V}_{x,[v]}=\{\xi\in T_{\tilde{X},(x,[v])}|\pi_*\xi\in \C\cdot v\}$$ for any $x\in X$ and any $v\in T_{X,x}\backslash \{0\}$.
Starting with a directed manifold $(X,V)=(X_0,V_0)$, we get a tower of directed manifolds $(X_k,V_k)$, called Demailly-Semple $k$-jet bundle of $X$, defined by $X_k=\tilde{X}_{k-1},V_k=\tilde{V}_{k-1}$. In particular, when $X$ is a hypersurface in $\CP^3$, we start with $(X, T_X)$

The line bundle $\ocal_{X_k}(1)$ will be called the Demailly-Semple $k-$jet line bundle.

\begin{theo}\cite{de}
The direct image sheaf $(\pi_{k,0})_*\ocal_{X_k}(m)$ on $X$ coincides
with the (locally free) sheaf $E_{k,m}V^*$ of $k$-jet differentials of weighted degree $m$, that
is, by definition, the set of germs of polynomial differential operators
\begin{equation} Q(f)=\sum_{\alpha_1\cdots\alpha_k}a_{\alpha_1\cdots\alpha_k}(f)(f')^{\alpha_1}(f'')^{\alpha_2}\cdots(f^{(k)})^{\alpha_k}
\end{equation}
on $J_kV$ (in multi-index notation, ($f')^{\alpha_1}=
((f_1')^{\alpha_{1,1}}(f_2')^{\alpha_{1,2}}\cdots(f_r')^{\alpha_{1,r}}) $, which are
moreover invariant under arbitrary changes of parametrization: a germ of operator
$Q\in  E_{k,m}V^*$ is characterized by the condition that, for every germ $f\in J_kV$ and
every germ $\phi \in \G_k$ ,
$$Q(f\circ\phi)
=
\phi'^mQ(f)\circ\phi$$
\end{theo}

Given a finit dimensional vector space $W$, we can define $E_{k,m}W^*$. Seeing that $W$ can be seen as the tangent space of its own, we define $E_{k,m}W^*$ to be the fibre at the origin of the bundle $E_{k,m}T_W^*$. More generally, if we have a vector bundle $V$ of finite rank over a manifold $X$. Let $G$ be the structure group of $V$, and let $\mathcal{P}$ be the principle bundle of $V$. We denote a fibre by $W$. An element $g\in G$ acts on $W$, inducing an automorphism of $E_{k,m}W^*$. Therefore, we get a representation of $G$ in $E_{k,m}W^*$. We can define $$\tilde{E}_{k,m}V^*\triangleq \mathcal{P}\times_GE_{k,m}W^*$$

However when $V$ is a subboudle of $T_X$, this vector bundle does NOT coincide with the original definition of $E_{k,m}V^*$. One can have a sense of this by taking a collection of covering charts of $X$ such that $V$ is trivial on each open set. Consider $V_1$ as a fiber bundle over $X$ with each fiber a vector bundle of rank $r$ over $\CP^{r-1}$. One can then compute the transition functions of $V_1$ on overlap of intersecting charts, and see that the second derivatives of the transition functions of $X$ are involved. The following example shows $E_{k,m}V^*$ is indeed not  a representation of the principle bundle of $V$.

Suppose $X$ is a surface, and $V=T_X$. In a neighborhood of a point $x\in X$ such that $T_X$ is trivial, one can easily compute that $E_{2,3}T_X^*|_x$ as a representation of $GL(T_{X,x})$ decomposes as a direct sum $S^3T_{X,x}^*\oplus K_{X,x}$. If $E_{2,3}T_X^*= \mathcal{P}\times_{GL{T_{X,x}}}E_{2,3}T_X^*|_x$, this decomposition will globalize to give a decomposition $E_{2,3}T_X^*=S^3T_X^*\oplus K_X$. But then we will have $H^0(X,E_{2,3}T_X^*\otimes (-K_X))\neq 0$

Recall that in \cite{dee} $\theta_{k,m}$ is defined to be the smallest rational number such that  $H^0(X,E_{k,m}T_X^*\otimes (tK_X))\neq 0$ , assuming $tK_X$ is an integral divisor, $t\in \Q$. Under our assumption, we get that $\theta_{2,3}\leq -\frac{1}{3}$. On the other hand, the following inequality was proved in \cite{dee}
\begin{prop}\cite{dee}
Let $X$ be a generic surface of degree $d\geq  6$ in $\CP^3$. Then
\begin{equation*}
\theta_{2,m}\geq \frac{-1}{2m}+\frac{2-7/(2m)}{d-4},\quad for \quad m=3,4,5
\end{equation*}

\end{prop}
We get a contradiction.

This example gives us the impression that the global structure of $E_{k,m}T_X^*$ is complicated. Therefore, knowing the local structure of $E_{k,m}T_X^*$ does not enable one to get global conclusions before one understand the transition functions for $E_{k,m}T_X^*$ on overlaps of open covering of $X$.

\end{document}